\documentclass[12pt]{article}  
\usepackage{latexsym,amsfonts}
\usepackage{amssymb}
\usepackage{amsbsy}
\usepackage{amsmath}
\usepackage{epsf}
\usepackage{graphicx}    
\newtheorem{theo} {Theorem}
\newcommand{\bteo}{\begin{theo}}  
\newcommand{\et}{\end{theo}}  
\newcommand{\bd}{\begin{displaymath}}  
\newcommand{\ed}{\end{displaymath}}  

\newcommand{\be} {\begin{equation}}  
\newcommand{\ee} {\end{equation}}  
\newcommand{\ba} {\begin{array}{l}}  
\newcommand{\ea} {\end{array}}  

\begin{document}

\begin{center}
 {\Large \bf
A reaction-diffusion  system \\ with nonconstant diffusion coefficients:\\
 exact and numerical  solutions}

\medskip

{\bf Roman Cherniha $^{a,b}$\footnote{\small  Corresponding author.
E-mail: r.m.cherniha@gmail.com; \\ roman.cherniha1@nottingham.ac.uk}
and Galyna Kriukova$^{b}$ }

$^{a}$ \quad School of Mathematical Sciences, University of Nottingham,\\
  University Park, Nottingham NG7 2RD, UK

$^{b}$  \quad National University of Kyiv-Mohyla Academy,\\ 2,
Skovoroda Street, Kyiv  04070, Ukraine.

\end{center}

\begin{abstract}

A  Lotka-Volterra type system with porous diffusion, which can be
used as an alternative model to the classical Lotka-Volterra system,
is under study. Multiparameter families of exact solutions of the
system in question
 are constructed and  their properties are established. It is shown that the solutions obtained can satisfy
  the zero Neumann conditions, which are  typical conditions for
  mathematical models describing real-world processes. It is proved
  that the system  possesses two stable steady-state  points
  provided  its  coefficients  are correctly-specified. In
  particular, this occurs when the system models  the prey-predator interaction.
  The exact solutions are used for solving boundary-value problems. The analytical  results are compared with  numerical solutions
  of the same   boundary-value problems but perturbed initial profiles.  It is  demonstrated  that
 the numerical  solutions  coincide
 with  the relevant
exact solutions with high exactness in the case of   sufficiently
small perturbations  of the initial profiles.

\end{abstract}


\section{ Introduction}

Nowadays  reaction-diffusion  systems are widely used in
mathematical modelling enormous variety of processes in ecology,
biology, medicine, physics, chemistry  and social sciences (see,
e.g., well-known books
  \cite{aris-75I, fife-79,ku-na-ei-16, britton,mur2002,
  mur2003,ch-dav-book,
okubo}). In the case of  interaction of  species (cells, chemicals,
etc.), the most popular reaction-diffusion  system is the diffusive
Lotka-Volterra system (DLVS) and its modifications. Extensive
studies of DLVS  started in the 1970s
 \cite{conway,hastings,carmi,rothe} and  new papers appear at
 regular bases (see \cite{ch-dav-2022} and references therein).
 Several generalisations of DLVS were worked out as well. The
 Shigesada--Kawasaki--Teramoto (SKT)
 model \cite{sh-ka-te}, involving non-constant
 coefficients of diffusion,
 is one of the most important. The recent rigorous results concerning
 this model can be found in \cite{pham-2019}--\cite{ch-da-ki-23}.

 The standard form of the SKT model reads as
 \be\label{0-1}
 \ba
  u_t = [(d_1+d_{11}u+d_{12}v)u]_{xx}+u(a_1-b_1u-c_1v),\\
 v_t = [(d_2+d_{21}u+d_{22}v)v]_{xx}+v(a_2-b_2u-c_2v),
 \ea
  \ee
where    $u=u(t,x)$  and $v=v(t,x)$ are two  unknown functions,
which usually  represent densities of two competing species (cells),
$d_1$ and $d_2$ denote the standard  diffusion coefficients,
$d_{11}u$ and $d_{22}v$ are intra-diffusion pressures,  $d_{12}v$
and $d_{21}u$ are cross-diffusion pressures, $a_1$ and $a_2$ are the
intrinsic growth coefficients, $b_1$ and $c_2$ denote the
coefficients of intra-specific competitions and $b_2$ and $c_1$
denote the coefficients of inter-specific competitions. Hereafter
the lower subscripts $t$ and $x$  denote differentiation with
respect to these variables. Notably, (\ref{0-1})  with $d_{ij}=0, \
i=1,2, \ j=1,2$ is nothing else but the well-known DLVS.

Neglecting cross-diffusion pressures (i.e. assuming that the
coefficients $d_{12}$ and $d_{21}$ are very small), keeping nonzero
$d_{11}$ and $d_{22}$ and replacing the nonlinear interaction  terms
$c_1uv$ and $b_2uv$ by the linear terms $h_1+c_1v$ and $h_2+c_2u$,
respectively, the SKT model is reduced to the form
\be\label{0-2}
 \ba
  u_t = [(d_1+d_{11}u)u]_{xx}+u(a_1-b_1u) -h_1-c_1v,\\
 v_t = [(d_2+d_{22}v)v]_{xx}+v(a_2-b_2v)-h_2-c_2u.
 \ea
  \ee
  Here new parameters $h_1$ and $h_2$ are introduced in order to
  take into account possible external forces (influences) on the
  interaction  of the species $u$  and $v$. Typically, such forces
  are ignored, i.e. $h_1=h_2=0.$
  So, the nonlinear system (\ref{0-2})  is a simplification of the
  classical SKT model, which is derived in a natural way under
  plausible assumptions. It can be noted that the substitution
 \be\label{0-3}
U= d_1+2d_{11}u, \ V= d_2+2d_{22}v \ee  transforms (\ref{0-2}) to
the equivalent  form
 \begin{equation}\ba
 U_t =(UU_x)_x+U(a^*_1+b^*_1U)+h^*_1+c^*_1V,\nonumber\\
 V_t =(VV_x)_x+V(a^*_2+b^*_2V)+h^*_2+c^*_2U,\ea
 \end{equation}
 where the parameters $a*_1, \ b*_1, \ h^*_1$ and $c^*_1$ are easily
 calculated via the parameters arising in (\ref{0-2}). In what
 follows the stars are skipped.
Thus, the  nonlinear reaction-diffusion (RD) system
\begin{equation}\label{3-1}\ba
 U_t =(UU_x)_x+U(a_1+b_1U)+h_1+c_1V,\\
 V_t =(VV_x)_x+V(a_2+b_2V)+h_2+c_2U, \ea
 \end{equation}
 will be the main object of this study. Obviously, we should assume
 that $c_1c_2\not=0$. Otherwise the system in question contains an
 autonomous equation, therefore its applicability would be  questionable.

  It should be noted this system
 can be regarded as a particular case  of the class of
extensively studied  RD systems with power nonlinearities (see,
e.g., \cite{ch-dav-book} and references therein).
For example, system (\ref{3-1})
 is a  Lotka-Volterra
type system with variable  diffusivities, in which the quadratic
terms are replaced by linear.
 On the other hand,  (\ref{3-1}) can be regarded as a
generalization of the porous-Fisher equation
\[ U_t =(UU_x)_x+U(1-U). \]
Physically, this equation is a model for the population dispersing
to regions of lower density more rapidly as the population gets more
crowded and one has been extensively studied  \cite [Section
13.4]{mur2002},\cite{witel-95}-- \cite{fadai-simp-2020}. So, the RD
system (\ref{3-1}) describes evolution of two populations with the
above habit that are additionally interacting according to a linear
low.

The  paper is organized as follows.
 In Section 2,  multiparameter families of exact solutions of the
 RD system (\ref{3-1})
 are constructed and  their properties are established.
 In the particular case, we show that the solutions obtained can satisfy
  the zero Neumann conditions, which are  typical conditions for
  mathematical models describing real-world processes.


 In Section 3, we find the stable  steady-state points  of  system (\ref{3-1})
  using the well-known procedures.
 We pay  main attention on the case when  the RD system (\ref{3-1}) possesses
 two stable positive nodes. It turns out, this case leads to very interesting behaviour of
 solutions, which we explore in the next section.
  Analysis  shows that the system coefficients must satisfy
 a very cumbersome system of nonlinear algebraic
inequalities if one aim to get two stable positive nodes. It is
demonstrated how to solve the inequalities obtained in the case of
the model describing the prey-predator interaction.

 In Section 4,  the analytical   results obtained in Sections 2-3
  are compared with numerical solutions obtained by simulations. The
  simulations were conducted using Python  scipy.integrate package.
The major conclusion of this section is the following: the  exact
solutions obtained  play an important role for solving some
boundary-value problems for  the  RD system (\ref{3-1}).  In fact,
the simulations show that
 the numerical  solutions  of boundary-value problems coincide
 with  the relevant
exact solutions with high exactness
in the case of   arbitrary  sufficiently small perturbations  of the
initial profiles generated by the exact solutions. The behaviour of
the numerical solutions in the case of  large  perturbations
   of the initial profiles are  studied as well. Finally, we present
some conclusions  in the last section.

\section{ Exact solutions of the RD system (\ref{3-1}) }


Plane wave solutions, in particular travelling waves,   are the most
common exact solutions,  which researcher are looking for.  Such
solutions for (\ref{3-1}) have  the form
\begin{equation} \label{2-2}
 U =\varphi(\omega), \quad  V =\psi(\omega),\quad \omega=x-kt \\
\end{equation}
where $k \in {\bf R}$ and the functions $\varphi$ and $\psi$ are
solutions of the  ODE system
\begin{equation}\label{2-3}\ba
(\varphi\varphi_\omega)_\omega+k\varphi_\omega+\varphi(a_1+b_1\varphi)+h_1+c_1\psi=0,\\
(\psi\psi_\omega)_\omega+k\psi_\omega+\psi(a_2+b_2\psi)+h_2+c_2\varphi=0.
\ea
\end{equation}  
The ODE system (\ref{2-3}) is not integrable  and only particular
solutions can be found.
The case $k=0$ of course  leads to time independent  solutions. To
the best of our knowledge, exact solutions of ODE system
(\ref{2-3}), in particular travelling fronts, are unknown. Some
exact solutions of the form (\ref{2-2}) follow from those
constructed in this section as particular cases.

 A much wider class of exact solutions of the RD system
(\ref{3-1}) can be constructed using the method of additional
generating conditions (MAGC) \cite{ch96}, \cite{ch98-c}, which is
related to the method of differential constraints
\cite{si-sha-yan-84,olv-94}. It can be noted that system (\ref{3-1})
is a particular case of a more general system, which was examined in
\cite{ch-king3} using MAGC. As additional generating conditions, the
 following   third order  ordinary differential equations (ODEs)
  \be\label{3-6}
 \ba
 \beta_1(t)
\frac{d u}{d x} + \beta_2(t) \frac{d^2 u }{ d x^2 }+
  \frac{d^3 u }{ d x^3}  =0\\
  \beta_1(t) \frac{d v }{ d x}+ \beta_2(t) \frac{d^2 v }{ d x^2 }+
  \frac{d^3 v }{ d x^3}  =0
\ea
 \ee 
were used. Here $ \beta_1(t)$ and $ \beta_2(t)$   are
to-be-determined smooth functions
 and the variable $t$ is considered as a parameter.
 It can be easily identified from  \cite{ch-king3} that the RD system
 (\ref{3-1}) with $b_1=b_2=b$ possesses exact solutions of the form
     \begin{equation}
U=\varphi_0(t)+\varphi_1(t)\exp(-\gamma x)+\varphi_2(t)\exp(\gamma
x)\\
 V=\psi_0(t)+\psi_1(t)\exp(-\gamma
x)+\psi_2(t)\exp(\gamma x)
 \label{2-6}
  \end{equation}
 if $b_1=b_2=-b<0$,
  and
    \begin{equation}
U=\varphi_0(t)+\varphi_1(t)\cos(\gamma x)+\varphi_2(t)\sin(\gamma
x)\\
 V=\psi_0(t)+\psi_1(t)\cos(\gamma x)+\psi_2(t)\sin(\gamma x),
 \label{2-7}
  \end{equation}
if if $b_1=b_2=b>0$, \  $\gamma= \sqrt{\frac{b}{2}}$. Notably,
ans\"atze  (\ref{2-6}) and (\ref{2-7})  follow from (\ref{3-6}) for
  $ \beta_1(t)=\pm\frac{b}{2}$ and $ \beta_2(t)=0$.

Direct calculations show that  ansatz  (\ref{2-6}) produces a family
of exact solutions for (\ref{3-1}) provided the functions
$\varphi_i$ and  $\psi_i \  (i=0,1,2)$ satisfy
  the ODE system
\begin{equation}\label{2-4}\ba
    \dot{\varphi_0}=-b\varphi_0^2
    +h_1+a_1\varphi_0+c_1\psi_0-2b\varphi_1\varphi_2,\\
    \dot{\psi_0}=-b\psi_0^2
    +h_2+a_2\psi_0+c_2\varphi_0-2b\psi_1\psi_2,\\
    \dot{\varphi_1}=-\frac{3}{2}b\varphi_0\varphi_1
    +a_1\varphi_1+c_1\psi_1,\\
    \dot{\psi_1}=-\frac{3}{2}b\psi_0\psi_1
    +a_2\psi_1+c_2\varphi_1,\\
    \dot{\varphi_2}= -\frac{3}{2}b\varphi_0\varphi_2
    +a_1\varphi_2+c_1\psi_2,\\
    \dot{\psi_2}=-\frac{3}{2}b\psi_0\psi_2
    +a_2\psi_2+c_2\varphi_2
\ea
    \end{equation}
    (here dots denote differentiation w.r.t. the time
    variable).
Similarly, ansatz  (\ref{2-7})  does the same provided the following
ODE system is satisfied
\begin{equation}\label{2-5}\ba
    \dot{\varphi_0}=b\varphi_0^2+h_1+a_1\varphi_0
    +c_1\psi_0+\frac{b}{2}(\varphi_1^2+\varphi_2^2),\\
    \dot{\psi_0}= b\psi_0^2+h_2+a_2\psi_0
    +c_2\varphi_0+\frac{b}{2}(\psi_1^2+\psi_2^2),\\
    \dot{\varphi_1}=\frac{3}{2}b\varphi_0\varphi_1
    +a_1\varphi_1+c_1\psi_1,\\
    \dot{\psi_1}=\frac{3}{2}b\psi_0\psi_1
    +a_2\psi_1+c_2\varphi_1,\\
    \dot{\varphi_2}= \frac{3}{2}b\varphi_0\varphi_2
    +a_1\varphi_2+c_1\psi_2,\\
    \dot{\psi_2}=\frac{3}{2}b\psi_0\psi_2
    +a_2\psi_2+c_2\varphi_2
\ea
    \end{equation}


Let us assume that interaction between two populations of species
takes place at the space interval $(0,L), \  L>0$ and the widely
used no-flux (zero Neumann ) conditions on boundary $L$
 \be\label{2-8}
  x=L: \quad  U_x=0, \quad V_x=0  \ee
  take place.
 Using  ansatz  (\ref{2-6}), one easily calculates that the no-flux
 conditions (\ref{2-8}) are satisfied
 iff
  \be\label{2-9}
   \varphi_2(t)=l\varphi_1(t), \quad \psi_2(t)=l\psi_1(t), \quad l=\exp(-2L\gamma). \ee
   Taking into account (\ref{2-9}),  the ODE system (\ref{2-4}) reduces to the form
 \begin{equation}\label{2-10}\ba
    \dot{\varphi_0}=-b\varphi_0^2
    +h_1+a_1\varphi_0+c_1\psi_0-2bl\varphi_1^2,\\
    \dot{\psi_0}=-b\psi_0^2
    +h_2+a_2\psi_0+c_2\varphi_0-2bl\psi_1^2,\\
    \dot{\varphi_1}=-\frac{3}{2}b\varphi_0\varphi_1
    +a_1\varphi_1+c_1\psi_1,\\
    \dot{\psi_1}=-\frac{3}{2}b\psi_0\psi_1
    +a_2\psi_1+c_2\varphi_1.
\ea
    \end{equation}
    So, we can formulate the following statement.
 \bteo
 The boundary-value problem (\ref{3-1}) with $b_1=b_2=-b<0$,
 (\ref{2-8}) and
 \be\label{2-8*}
  x=0: \quad  U =\varphi_0(t)+(1+l)\varphi_1(t), \quad V =\psi_0(t)+(1+l)\psi_1(t).  \ee
posses the exact solution
  \begin{equation}\label{2-11}
 \ba
U=\varphi_0(t)+\varphi_1(t)(\exp(-\gamma x)+l\exp(\gamma
x))\\
 V=\psi_0(t)+\psi_1(t)(\exp(-\gamma
x)+l\exp(\gamma x)), \ l=\exp(-2L\gamma) \ea
  \end{equation}
  provided the functions  $\varphi_k, \psi_k, k=1,2$
  form a  solution of the ODE system (\ref{2-10}).
  \et
 {\bf Remark 1.} This result is also correct  if one
 sets  $L= +\infty$, i.e. $ l=0 $, however,   cannot be extended on the case when  solution
 (\ref{2-11}) satisfies the zero Neumann conditions on both boundaries
 of the interval $(0,L)$.

It can be noted that system (\ref{2-10}) with $l=0 $  reduces to the
form

\begin{equation}\label{4-5} \ba
\dot{\varphi}_0=-b\varphi_0^2+h_1+a_1\varphi_0+c_1\psi_0\\
\dot{\psi}_0=-b\psi_0^2+h_2+a_2\psi_0+c_2\varphi_0\\
\dot{\varphi}_1=-\frac
32b\varphi_0\varphi_1+a_1\varphi_1+c_1\psi_1\\
\dot{\psi}_1=-\frac 32b\psi_0\psi_1+a_2\psi_1+c_2\varphi_1.
\ea  \end{equation}  

Assuming  that  $(U_0,V_0)$ is steady-state solution of the RD
system (\ref{3-1}) with  $b_1=b_2=-b<0$, i.e

\begin{equation}\label{4-6}
\ba   b U_0^2=h_1+a_1U_0+c_1V_0\\ b V_0^2=h_2+a_2V_0+c_2U_0, \ea
 \end{equation}  
 the nonlinear ODE system (\ref{4-5}) with $(\varphi_0,\psi_0)=(U_0,V_0)$ reduces to
 the linear  system
\begin{equation}\label{4-7}
\ba  \dot{\varphi}_1=\Big(a_1-\frac
32bU_0\Big)\varphi_1+c_1\psi_1\\
\dot{\psi}_1=c_2\varphi_1+\Big(a_2-\frac 32 bV_0\Big)\psi_1.
 \ea
   \end{equation}  

  So we obtain  the following general solutions of
(\ref{4-7}):

if $\triangle =[(a_1-a_2)+\frac 32b(V_0-U_0)]^2+4c_1c_2>0$ then

\begin{equation} \label{4-8}
\ba   \varphi_1=e_1c_1\exp(s_1t)+e_2\Big(s_2-a_2+\frac 32
bV_0\Big)\exp(s_2t)\\ \psi_1=e_1\Big(s_1-a_1+\frac 32
bU_0\Big)\exp(s_1t)+e_2c_2\exp(s_2t) \ea
 \end{equation}  
where $s_{1,2}=\frac 12\Big[a_1+a_2-\frac 32
b(U_0+V_0)\pm\sqrt{\triangle}\Big]$;

if $\triangle<0$ then

\begin{equation}\label{4-9}
 \ba
\varphi_1=c_1\exp(pt)(e_1\sin(q t)+e_2 \cos(q t))\\
\psi_1=\exp(pt)\left[(e_1(p-a_1+\frac 32 bU_0)-e_2q)\sin(q
t)+\right.\\ \left.\qquad+(e_1q+e_2(p-a_1+\frac 32 bU_0))\cos(q
t)\right]
 \ea
  \end{equation}  
  where $p=\frac 12\Big(a_1+a_2-\frac 32 b(U_0+V_0)\Big),\quad q=\frac
  12\sqrt{-\triangle}$;

   if $\triangle=0$ then
\begin{equation}\label{4-10}
 \ba
\varphi_1=(e_1c_1t+c_2)\exp(st)\\
\psi_1=\left[e_1(1\pm\sqrt{-c_1c_2}t)\pm
e_2\sqrt{\frac{-c_2}{c_1}}\right]\exp(st)
 \ea
  \end{equation}  
  where $s=\frac 12\Big[a_1+a_2-\frac 32 b(U_0+V_0)\Big]$. Here $e_1$ and
  $e_2$ are arbitrary constants.

Thus, three  families of exact solutions of the form

\begin{equation}\label{4-14}
 \ba
U=U_0+\varphi_1(t)\exp\Big(-\sqrt{\frac b2}x\Big)\\
V=V_0+\psi_1(t)\exp\Big(-\sqrt{\frac b2}x\Big)
 \ea
\end{equation} for  the RD system
\begin{equation}\label{2-12}
\ba
 U_t =(UU_x)_x+U(a_1-b U)+h_1+c_1V\\
 V_t =(VV_x)_x+V(a_2-b V)+h_2+c_2U
 \ea
  \end{equation}
 are  constructed.


One observes that all solutions of the form (\ref{4-14}) satisfy the
zero Neumann boundary conditions (\ref{2-8}) with $L= +\infty$. The
asymptotical behaviour   (with respect to the time) of these
solutions depends essentially on parameters  $s_1, s_2 $ (see
(\ref{4-8})),  $p$ (see (\ref{4-9})) and $s $ (see (\ref{4-10})). If
these parameters are negative then the relevant solution
(\ref{4-14}) tends to the steady-state solution $(U_0,V_0)$ for $t
\to +\infty$, otherwise one  is an unbounded solution   or  a
periodical solution (see (\ref{4-9})with $p=0$).


Now we consider the ansatz  (\ref{2-7}). This ansatz can be examined
in quite similar way, however, the result is  different. So, we
consider again  interaction between two populations of species and
take the  interval $(L_1, L_2)$, where $
L_1=L+\frac{k_1\pi}{\gamma}, \, L_2 =L+\frac{k_2\pi}{\gamma},
 \, k_1< k_2 \in {\bf Z}$. We may set
 $l=\tan(L_1\gamma)=\tan(L_2\gamma)$  and  no-flux  conditions on the both
 boundaries
 \be\label{2-8**}
  x=L_1,\  x=L_2: \  U_x=0, \ V_x=0.  \ee

  \bteo
 The boundary-value problem (\ref{3-1}) with $b_1=b_2=b>0$ and
 (\ref{2-8**})
posses the exact solution
  \begin{equation}\label{2-13}
 \ba
U=\varphi_0(t)+\varphi_1(t)(\cos(\gamma x)+l\sin(\gamma x))\\
 V=\psi_0(t)+\psi_1(t)(\cos(\gamma x)+l\sin(\gamma x)),
\ea
  \end{equation}
   provided the functions  $\varphi_k, \psi_k, k=1,2$
  form a  solution of the ODE system
  \begin{equation}\label{2-14} \ba
    \dot{\varphi_0}=b\varphi_0^2+h_1+a_1\varphi_0
    +c_1\psi_0+\frac{b}{2}(1+l^2)\varphi_1^2,\\
    \dot{\psi_0}= b\psi_0^2+h_2+a_2\psi_0
    +c_2\varphi_0+\frac{b}{2}(1+l^2)\psi_1^2, \\
    \dot{\varphi_1}=\frac{3}{2}b\varphi_0\varphi_1
    +a_1\varphi_1+c_1\psi_1,\\
    \dot{\psi_1}=\frac{3}{2}b\psi_0\psi_1
    +a_2\psi_1+c_2\varphi_1.
\ea
    \end{equation}
  \et

 It is very unlikely that  the nonlinear ODE system (\ref{2-14}) is  integrable for
  arbitrary coefficients, however,
 we were able to find  its particular solutions,  setting $l=1$
 (this restriction is only for convenience)
     and $ \varphi_0-\psi_0 = \textrm{const},  \varphi_1=-\psi_1$.
As a result, the following exact solutions
  were constructed:

\be\label{2-15}
    \ba
U= \frac{2}{3b(t_0-t)}(1\pm\cos(\sqrt{\frac b2}(x-x_0))-\frac{4c_1}{b} \\
 V=   \frac{2}{3b(t_0-t)}(1\mp\cos(\sqrt{\frac b2}(x-x_0))-\frac{4c_2}{b}
 \ea \ee
 and
 \be\label{2-16}
    \ba
U=\frac{a_1-7c_1}{3b}\Bigl(\tanh\Bigl(\frac{a_1-7c_1}{2}(t_0-t)\Bigr)
-1\Bigr)(1\pm\cos(\sqrt{\frac b2}(x-x_0))-\frac{4c_1}{b},\\
 V=  \frac{a_1-7c_1}{3b}\Bigl(\tanh\Bigl(\frac{a_1-7c_1}{2}(t_0-t)\Bigr)
-1\Bigr)(1\mp\cos(\sqrt{\frac b2}(x-x_0))-\frac{4c_2}{b}
 \ea \ee
 of the nonlinear RD systems
   \be\label{2-17}
    \ba
U_t=(U U_x)_x+U(7c_1+b U)+c_1V+\frac{4c_1}{b}(3c_1+c_2)\\
V_t=(V V_x)_x+V(7c_2+b V)+c_2U+\frac{4c_2}{b}(3c_2+c_1)
\ea
\ee
and
   \be\label{2-18}
    \ba
U_t=(U U_x)_x+U(a_1+b U)+c_1V+\frac{4c_1}{b}(a_1-4c_1+c_2) \\
V_t=(V V_x)_x+V(a_2+b V)+c_2U+\frac{4c_2}{b}(a_2-4c_2+c_1),
    \ea
\ee respectively. Here  $t_0$ and $x_0$ are arbitrary constants, $
a_1-7c_1=a_2-7c_2 \not=0$.

Note that solution (\ref{2-16}) is also valid if one replaces the
function $\tanh$ by $\coth$.
Interestingly,  solutions (\ref{2-15}) and (\ref{2-16}) with the
function $\coth$ instead of $\tanh$ blow-up  for a finite time
$t_0>0$. This is an essential difference from the  solutions
obtained for the RD system (\ref{3-1}) with $b_1=b_2=-b<0$.



\section{ The spatially-homogeneous  case }

 In this section,  constant steady-state points of the RD system
 (\ref{3-1}) with $b_1=b_2=-b<0$ are  studied.
Since we are looking for constant steady-state points, the  ODE
system
 \begin{equation}\label{5-2}\ba
 U_t = U(a_1  - b U)+h_1+ c_1V,\\
 V_t =  V(a_2  - b V) +h_2+c_2U.
\ea\end{equation}
 should be considered instead of the  RD system
 (\ref{3-1}).
   Any steady state  point of (\ref{5-2})
($U_i, V_i$) must  satisfy the system of algebraic equations
 \begin{equation}\label{5-3}\ba
U_i(a_1  - b U_i)+h_1+ c_1V_i = 0,\\
V_i(a_2  - b V_i) +h_2+ c_2U_i = 0. \ea\end{equation}

 It is well-known  that the character of
a nonsingular steady-state point can be established by linearization
of the right-hand-side (RHS) of (\ref{5-2}). In fact, taking into
account only linear terms of Taylor's series in the point ($U_i,
V_i$), we receive the algebraic system
\[\ba
(a_1  - 2 b U_i)(U - U_i)+c_1(V - V_i) = 0, \\  c_2(U - U_i)+(a_2 -
2 b V_i)(V - V_i) = 0.\ea
\]

 The     eigenvalues $\lambda_1$ and $\lambda_2$  of the matrix
\begin{equation}\label{5-5}
 [M]= \left(%
\begin{array}{cc}
  a_1-2bU_i & c_1 \\
  c_2 & a_2-2bV_i \\
\end{array}%
\right).
\end{equation}
can be calculated by formulas
\begin{equation}\label{5-7} \lambda_{1,2} =  \frac{1}{2}(A \pm
\sqrt{A^2-4B)},  \ A=a_1+a_2 - 2b(U_i + V_i),\
B=(a_1-2bU_i)(a_2-2bV_i) -c_1c_2.
\end{equation}

The type of the steady-state point $(U_i, V_i)$ is determined by the
sign of the $ Re\lambda _1$ and $ Re\lambda _2$. Depending on $(U_i,
V_i)$  and the coefficients arising in (\ref{5-2}), one can expect
to get a wide range of types of steady-state points. Obviously, the
algebraic system (\ref{5-3}) has up to four real solutions, however,
only those with non-negative coordinates $U_i$ and  $V_i$  are
interesting from applicability point of view.  In the case of the
classical Lotka-Volterra competition model and the SKT model, the
most interesting phenomenon occurs when there is a stable
steady-state points with two non-zero coordinates (see \cite{lou-ni}
and references cited therein). This case is interpretable as
coexistence of two population of species (cells).

Here it is shown that, in contrast to the above mentioned models,
the  ODE system (\ref{5-2}) possesses  {\it two stable steady-state
points with positive coordinates} provided its coefficients are
correctly defined. For such  stable steady-state points,
 the following inequalities should take place
\begin{equation}
Re\lambda _1<0, \quad Re\lambda _2<0, \quad  U_i >  0, \quad V_i
>  0, \quad i=1,2.
\label{5-10}\end{equation} In order to avoid cumbersome
calculations, we assume  the additional restrictions
\begin{equation}
V_i=\alpha U_i, \quad  \quad i=1,2 \label{5-10a}
\end{equation} where $\alpha >0$ is a given number.
Using  (\ref{5-3}), one obtains
\[
b U_i^2-(a_1+c_1\alpha )U_i -h_1= 0,
\] \[
 b U_i^2- {{a_2\alpha
+c_2} \over {\alpha ^2}}U_i - {{h_2} \over {\alpha ^2}}  = 0,
\]
hence
\begin{equation}\label{5-13}\ba
 b U_i^2-(a_1+c_1\alpha )U_i -h_1= 0, \\ a_1+c_1\alpha =
\frac{a_2\alpha +c_2}  {\alpha ^2}, \ h_2=h_1\alpha ^2.
\ea\end{equation}

The  first equation of (\ref{5-13}) gives
\begin{equation}
\label {5-16} U_{1,2} = \frac{1}{2b}(a_1+c_1\alpha \pm
\sqrt{(a_1+c_1\alpha)^2+4bh_1}), \end{equation}
therefore  we arrive at  the restrictions
\begin{equation}
a_1+c_1\alpha > 2 \sqrt{b|h_1|}, \quad h_1 <0, \quad h_2<0
\label{5-17}\end{equation}
 in order to get  only positive values of
$U_{1,2}$. Substituting   (\ref{5-16}) into (\ref{5-7}), we obtain
\[
A_\pm  = -(c_1\alpha +\frac{c_2}{\alpha }) \mp(1+\alpha
)\sqrt{(a_1+c_1\alpha )^2+4bh_1}
\]
and
\[
B_\pm  = \alpha [(a_1+c_1\alpha )^2+4bh_1]\pm (c_1\alpha^2
+\frac{c_2}{ \alpha }) \sqrt{(a_1+c_1\alpha )^2+4bh_1}.
\]

Thus, taking into account (\ref{5-10}), (\ref{5-16}) and
(\ref{5-17}),   we need to satisfy the following system of
inequalities:

\begin{equation}\label{5-20}\ba
{1\over\alpha }|c_1\alpha ^2 + \frac{c_2}{ \alpha }| < \sqrt{
(a_1+c_1\alpha)^2+4bh_1} <\frac{1}{1+\alpha} (c_1 \alpha +{c_2\over
\alpha}),
 \\ a_1+c_1\alpha > 2\sqrt{b|h_1|}, \quad b>0, h_1 <0,
h_2 <0, \ea \end{equation} where
\[
\alpha=\sqrt{\frac{h_2}{h_1}}, \quad a_1+c_1\alpha = a_2\alpha
+\frac{c_2}{ \alpha^2 }.\]

 One can easily  derive from inequalities
 (\ref{5-20}) that the
coefficients $c_1 $ and $c_2$ must satisfy the inequality
\begin{equation}
c_1c_2<0.
 \label{5-22}\end{equation}
 Thus,  one needs  to
consider  two different cases: {\it (i)} $ c_1<0, c_2>0$ and
{\it(ii)} $c_1>0, c_2<0$. In the first case, inequalities
(\ref{5-20}) lead to the requirement $ \alpha >1 $, while $0< \alpha
<1 $ is obtained in the second that. Note that the special case $
\alpha =1 $ leads to contradiction (see the first line in
(\ref{5-20})).  Now one notes that case {\it(ii)} can be reduced to
{\it(i)} by the following renaming
\[ U \to V, \, V \to U, \, a_i \to a_{3-i}, \, b_i \to b_{3-i}, \, c_i \to c_{3-i}, \,
h_i \to h_{3-i}, \, i=1,2 \] in the ODE system  (\ref{5-2}). So it
is enough to consider only case {\it (i)}.

It is well-known that the precise type of a stable point depends on
the      eigenvalues $\lambda_1$ and
 $\lambda_2$. If both lambda-s are real then stable nodes are
 obtained, while stable spires arise for complex lambda-s.

  To get  two positive
stable {\it nodes}, we need  to  satisfy  the  conditions
\begin{equation}
 A_+^2 >4B_+, \quad  A_-^2 >4B_-.
\label{5-26}\end{equation}
 Conditions (\ref{5-26}) are equivalent to the inequality
\begin{equation}
\Lambda _2s^2-2\Lambda _1s+\Lambda _0>0, \label{5-27}\end{equation}
where
\[
 \Lambda _0=(c_1 \alpha +\frac{c_2}{\alpha})^2, \
\Lambda _1=|c_1 \alpha^2+\frac{c_2}{\alpha}-c_1\alpha -c_2|,
 \]
 \[
\Lambda _2=(1-\alpha )^2, \  s=\sqrt{(a_1+c_1\alpha )^2+4bh_1}.
 \]

The solutions of  inequality (\ref{5-27}) are determined by the
roots
 \[s_{1,2}={c_1\alpha -{c_2\over\ \alpha}\pm
2\sqrt{-c_1c_2} \over\ \alpha -1 } ={(\sqrt{-c_1\alpha } \pm
\sqrt{{c_2\over\ \alpha}})^2 \over\ \alpha -1 } \] of equation
$\Lambda _2s^2-2\Lambda _1s+\Lambda _0=0$.  So,  we obtain the
restrictions
\begin{equation}
s<{(\sqrt{-c_1\alpha }-\sqrt{{c_2\over\ \alpha}})^2 \over\ \alpha
-1 } \qquad \mbox{  or}\; \qquad \mbox{ }\;
 s>{(\sqrt{-c_1\alpha }+\sqrt{{c_2\over\ \alpha}})^2 \over\ \alpha -1 } .
\label{5-31}\end{equation}
 To
obtain two positive stable nodes, one needs to solve
 system (\ref{5-20}) together with   inequalities
(\ref{5-31}). Taking into account (\ref{5-20}), it is easily seen
that the second inequality in (\ref{5-31}) cannot be satisfied
because
\[ \frac{(\sqrt{-c_1\alpha }+\sqrt{\frac{c_2}{\alpha}})^2}
{\alpha -1 }> \frac{1}{1+\alpha } (c_1 \alpha +\frac{c_2}{\alpha}).
\] Examining the first  inequality in (\ref{5-31}), we arrive at the
system of algebraic restrictions
\begin{equation} \label{5-34}
\sqrt{4b|h_1|+ (c_1\alpha  + \frac{c_2}{\alpha^2})^2}-c_1\alpha <
a_1< \sqrt{4b|h_1|+\beta^2}-c_1\alpha,
\end{equation} where
\begin{equation}\label{5-34a}\ba
 b>0, \quad  h_1 <0,  \quad h_2 <0,
\quad c_1<0,  \quad c_2>0
  \\ \beta = \textrm{min}\Big(\frac{1}{1+\alpha } |c_1
\alpha +\frac{c_2}{\alpha}|, \frac{(\sqrt{-c_1\alpha
}-\sqrt{\frac{c_2}{\alpha}})^2 }{ \alpha -1 } \Big), \\ \alpha
=\sqrt{\frac{h_2}{h_1}}>1, \, a_1+c_1\alpha = \frac{a_2\alpha
+c_2}{\alpha^2}. \ea
 \end{equation}

{\it Thus, ODE system (\ref{5-2}) possesses  two stable  nodes with
the positive coordinates (\ref{5-10a}) provided its coefficients $b,
a_i, c_i$ and $h_i, i=1,2$,
   satisfy the algebraic inequalities  (\ref{5-34})--(\ref{5-34a}).}

   To get  two positive
stable {\it spires}, we need  to  analyze  the  inequalities
\begin{equation}
 A_+^2 <4B_+, \quad  A_-^2 <4B_-,
\label{5-35}\end{equation} which  can be solved in a quite similar
way. As a result, one obtains
\begin{equation} \label{5-36}
\sqrt{4b|h_1|+ \gamma^2}-c_1\alpha < a_1<
\sqrt{4b|h_1|+\frac{1}{(1+\alpha)^2 } (c_1 \alpha
+\frac{c_2}{\alpha})^2}-c_1\alpha,
\end{equation} where
\[ \gamma= \max\Big(|c_1
\alpha +\frac{c_2}{\alpha^2}|, \frac{(\sqrt{-c_1\alpha
}-\sqrt{\frac{c_2}{\alpha}})^2 }{ \alpha -1 } \Big), \] while other
restrictions are the same as in (\ref{5-34a}).


It should be noted  that assumptions (\ref{5-10a}) do not  allow us
to get the  ODE system (\ref{5-2}) with $h_1=h_2=0$  possessing two
stable steady-state points with positive coordinates. In fact,
according to (\ref{5-17})
 one needs the restriction $h_1h_2 \not=0.$
 However, it can be shown that the  ODE system (\ref{5-2}) with $h_1=h_2=0$
 possesses such steady-state points if one skip (\ref{5-10a}).
 In fact, we may write
 \begin{equation}
V_i=\alpha_i U_i, \ \alpha_i>0,  \quad i=1,2 \label{5-37}
\end{equation}
in the general case. Solving (\ref{5-3}) with $h_1=h_2=0$
(straightforward calculations are omitted here), one obtains two
steady-state points with the positive  coordinates
\begin{equation}
U_i=\frac{1}{b}(a_1+ \frac{c_1}{\alpha_i}),  \quad   V_i=
\frac{1}{\alpha_i b}(a_1+ \frac{c_1}{\alpha_i}), \  i=1,2, \quad
a_1>\max\{ \frac{|c_1|}{\alpha_1},\frac{|c_1|}{\alpha_2}\}
\label{5-39*}
\end{equation}
provided  $\alpha_i$ are positive roots of the polynomial $P_3(y)=
c_2y^3+a_2y^2-a_1y-c_1$. A simple analysis shows that two positive
roots exist, for example, if the following restrictions take place:
\begin{equation} a_1>0, \ c_1<0, \ c_2>0. \label{5-38*}
\end{equation} The  ODE system (\ref{5-2}) with
$h_1=h_2=0$ and the above restrictions can be considered as a model
for prey-predator interaction.

Now one needs to satisfy (\ref{5-10}) in order to obtain two stable
steady-state points. Actually, the last two inequalities are already
guaranteed, while for the first two, one needs
\begin{equation}\label{5-39} \ba
\lambda_{1,2} =  \frac{1}{2}(A \pm \sqrt{A^2-4B}),  \\
 A=a_2-a_1 -  \frac{2}{\alpha_i}\Big(a_1+c_1(1+
 \frac{1}{\alpha_i})\Big)<0, \\
 B=(a_1+ \frac{2c_1}{\alpha_i})\Big(\frac{2}{\alpha_i}(a_1+ \frac{c_1}{\alpha_i})-a_2\Big)-c_1c_2>0. \ea\end{equation}
The inequalities in (\ref{5-39}) are equivalent to the coefficient
restrictions
\begin{equation}\label{5-40} \ba
a_2<  2(1+ \frac{1}{\alpha_i})(a_1+ \frac{c_1}{\alpha_i})- a_1\\
 c_2> (\frac{a_1}{c_1}+\frac{2}{\alpha_i})\Big(\frac{2}{\alpha_i}(a_1+\frac{c_1}{\alpha_i})-a_2\Big).\ea\end{equation}

 {\it Thus, the ODE system (\ref{5-2}) with $h_1=h_2=0$  possesses  two stable positive
steady-state points (\ref{5-39*})  provided its coefficients
   satisfy the algebraic inequalities  (\ref{5-40}), where $\alpha_i \
   i=1,2$ are positive roots of $P_3(y)=
c_2y^3+a_2y^2-a_1y-c_1$. The latter is guaranteed by the
restrictions (\ref{5-38*}).}

A simple example of such system occurs if one assumes that
$\alpha_1=1, \ \alpha_2=\frac{1}{c_2}>0$. Indeed, simple
calculations show that the ODE system
\begin{equation}\label{5-2*}\ba
 U_t = U\Big((4+\frac{5}{c_2})  - b U\Big)-\frac{5}{c_2}V,\\
 V_t =  V\Big((4-c_2)  - b V\Big) +c_2U
\ea\end{equation} with $b>0$ and $0<c_2 \leq \frac{1}{2}$  possesses
two stable positive steady-state points
$P_1=(\frac{4}{b},\frac{4}{b})$ and
$P_2=(\frac{5-c_2}{bc_2},\frac{5-c_2}{b}).$


In conclusion of this section,   it should be noted that all results
obtained above
 are valid  also for any RD systems of the form
 \begin{equation}\label{5-1} \ba
 U_t = (D_1(U,V)U_{x})_x +U(a_1  - b U)+h_1+ c_1V,\\
 V_t =(D_2(U,V)V_{x})_x  +V(a_2  - b V) +h_2+ c_2U,
\ea\end{equation} where   $D_1$ and $D_2$ are variable diffusivities
given by smooth nonnegative functions. In the case of  constant
diffusivities $D_1$ and $D_2$, one can set $h_1=h_2=0$ in
(\ref{5-1}) without losing generality. In fact, the simple
substitution
\[ U=U^*+d_1, \quad  V=V^*+d_2, \]
where the constants $d_1$ and $d_2$ are the solutions of the
algebraic system
\[ bd_1^2=a_1d_1+c_1d_2+h_1, \quad
bd_2^2=a_2d_2+c_2d_1+h_2, \] reduces (\ref{5-1}) to the same system
for $U^*$ and $V^*$ with $a_k^*=a_k-2bd_k, \, c_k^*=c_k, \, b^*=b$
and  $h_1^*=h_2^*=0$. However, this substitution does not work in
the case of non-constant $D_1$ and $D_2$,  hence  one cannot  drop
$h_1$ and $h_2$ in the RD system (\ref{5-1}) without losing
generality.

\section{ Properties of solutions  of a RD system with two stable
steady-state points}


 In this section, we  investigate the  properties  of solutions
 of   the RD system   (\ref{3-1})
 with
  $b_1=b_2=-b<0$, assuming that the relevant dynamical system
  possesses two stable nodes.
 First of all, we need to    construct an example of  such system
 with correctly-specified coefficients using the results of Section
 3.
  Let us choose $\alpha =3, h_1=-1 $ and $c_1=-1$, therefore
(\ref{5-34a}) immediately gives  $h_2=-9$.
 Setting now $c_2=27$, we observe that $\xi=0$ and $\eta=
 3(2-\sqrt{3})$ (see again (\ref{5-34a})).
 Setting $b=3$,  one obtains $3+2\sqrt{3} < a_1 < 3+\sqrt{75-36\sqrt{3}}$.
So, choosing, for example,
 $a_1=6.5$, we  obtain $a_2$=1.5 (see the last   equation of
(\ref{5-34a})).   Thus,  two  stable  nodes are
 $P_1=(U_1, V_1) = (\frac{1}{2},\frac{3}{2})$ and $P_2=(U_2, V_2)
=(\frac{2}{3},2)$ (see (\ref{5-10a}) and (\ref{5-16})). Using
(\ref{5-3}) and $P_i,\ i=1,2$,  the third and fourth steady-state
points: $P_3=(U_3, V_3) = ({19-\sqrt{145}\over 12},
(\frac{\sqrt{145}-5)}{4})\approx(0.58, 1.76)$ and $P_4=(U_4, V_4)
=(\frac{19+\sqrt{145}}{12}, \frac{-5-\sqrt{145}}{4})\approx (2.59,
-4.26)$ were determined.
 Calculating the eigenvalues
 $\lambda_1$ and $\lambda_2$ of matrix (\ref{5-5})
  for the  points $(U_3, V_3)$ and $(U_4, V_4)$,
 we obtain   $\lambda_1<0, \lambda_2>0$ in both cases, therefore
$(U_3, V_3)$ and $(U_4, V_4)$ are the saddle points.

Thus, the RD system (\ref{3-1}) with the above specified
coefficients reads as
 \begin{equation}\label{5-36}\ba
 U_t = (UU_{x})_x +U(6.5  - 3U)-1-V \\
 V_t = (VV_{x})_x  +V(1.5  - 3V) -9+ 27U.\ea
\end{equation} From mathematical point of view,
(\ref{5-36}) is a system of two
 porous-Fisher equations with additional  linear reaction  terms.
 This system can be considered as a model for prey-predator
 interaction of the species $U$ and $V$ (these  functions represent non-dimensional densities
  of preys and
 predators, respectively). The constants $-1$ and $-9$  can be
 thought as a removal of fixed numbers of preys and
 predators
  by an external force, e.g., by human.

\begin{figure}[t]
\begin{minipage}[t]{12.0cm}
\centerline{\includegraphics[width=12.0cm]{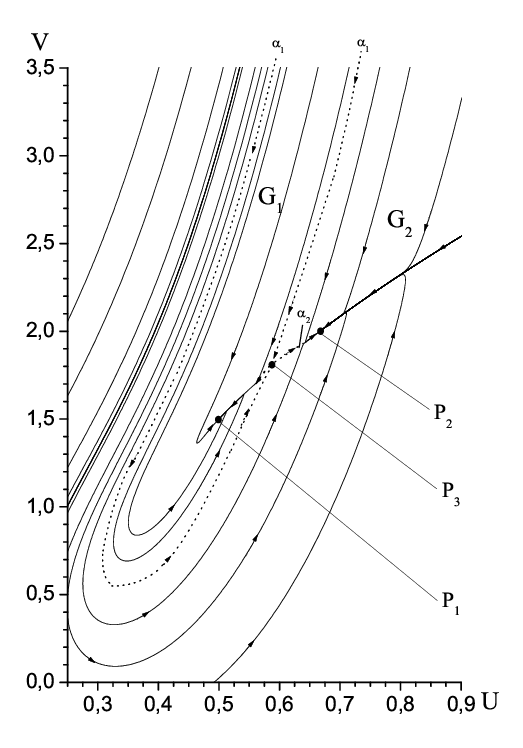}} \caption{The
first quadrant of the  $(U,V)$-phase plane of solutions for the RD
system (\ref{5-36}) with $(UU_{x})_x =(VV_{x})_x =0$.}
\label{Fig-2b}
\end{minipage}
 \end{figure}

Now we consider  the  dynamical   system,
generated by the system (\ref{5-36})  with  $(UU_{x})_x =(VV_{x})_x
=0$.
 Nowadays the relevant  $(U,V)$-phase
plane of solutions can easily  be constructed  using many  existing
program packages.  As a result, we derived  $(U,V)$-phase plane
presented in  Fig.\ref{Fig-2b}. Of course, the first quadrant only
is interesting from applicability point of view because the
functions $U$ and $V$ must be nonnegative. Three of four
steady-state points satisfy this requirement, while the unstable
steady-state point $P_4$ belongs to another quadrant and has no
biological meaning.

  It can be noted that the
 separatrix $\alpha_1$(see the dotted line) divides  the first
quadrant into two  domains $G_1$ and $G_2$. The  domain $G_1$
contains only curves leading to the stable node $P_1$, while $G_2$
contains only those leading to  $P_2$. The saddle point $P_3$ is a
cross-point of the  separatrix $\alpha_1$ and the separatrix
$\alpha_2$ (the second one is only partly pictured in
Fig.\ref{Fig-2b}).

Now we turn to exact solutions of the  RD system (\ref{5-36}).
Having the known  steady-state points and using formulas
(\ref{4-8})-(\ref{4-10}) and (\ref{4-14}), we can construct four
families exact solutions in the explicit form. Thus, the
 stable  nodes  $P_1=(U_1, V_1)$ and  $P_2=(U_2, V_2)$ generate the
 two-parameter  families exact solutions
\begin{equation}\label{5-37}
 \ba
U={1\over 2}-\exp\Big(-\sqrt{\frac {3}{2}}x-{1\over 2}t
\Big)\biggl(e_1\sin(\frac {\sqrt{71}}{4}t)+e_2 \cos(\frac
{\sqrt{71}}{4}t)\biggr)\\
 V={3\over 2}- {1\over 4}\exp\Big(-\sqrt{\frac {3}{2}}x-{1\over 2}t
\Big)\biggl(\Big(19e_1+\sqrt{71}e_2\Big)\sin(\frac {\sqrt{71}}{4}t)\\
\qquad+\Big(19e_2-\sqrt{71}e_1\Big) \cos(\frac
{\sqrt{71}}{4}t)\biggr)
 \ea
\end{equation} and
\begin{equation}\label{5-38}
 \ba
U={2\over 3}-e_1\exp\Big(-\sqrt{\frac {3}{2}}x-\frac
{8-\sqrt{13}}{4}t \Big)+\frac {22-\sqrt{13}}{4}
e_2\exp\Big(-\sqrt{\frac {3}{2}}x-\frac {8+\sqrt{13}}{4}t \Big)
\\ V=2- \frac {22-\sqrt{13}}{4}e_1\exp\Big(-\sqrt{\frac
{3}{2}}x-\frac {8-\sqrt{13}}{4}t \Big)+ 27e_2\exp\Big(-\sqrt{\frac
{3}{2}}x-\frac {8+\sqrt{13}}{4}t \Big), \ea
\end{equation} for $P_1$ and $P_2$,  respectively. In
fact, $\triangle$ is negative for $P_1$, while $\triangle>0$ for
$P_2$, as a result, the exact solutions (\ref{5-37}) and
(\ref{5-38}) have essential different structures.  Here $e_1$ and
  $e_2$ are arbitrary constants, which can be specified according
  to the given initial profiles.

\begin{figure}[t]
\begin{minipage}[t]{9.0cm}
\centerline{\includegraphics[width=9.0cm]{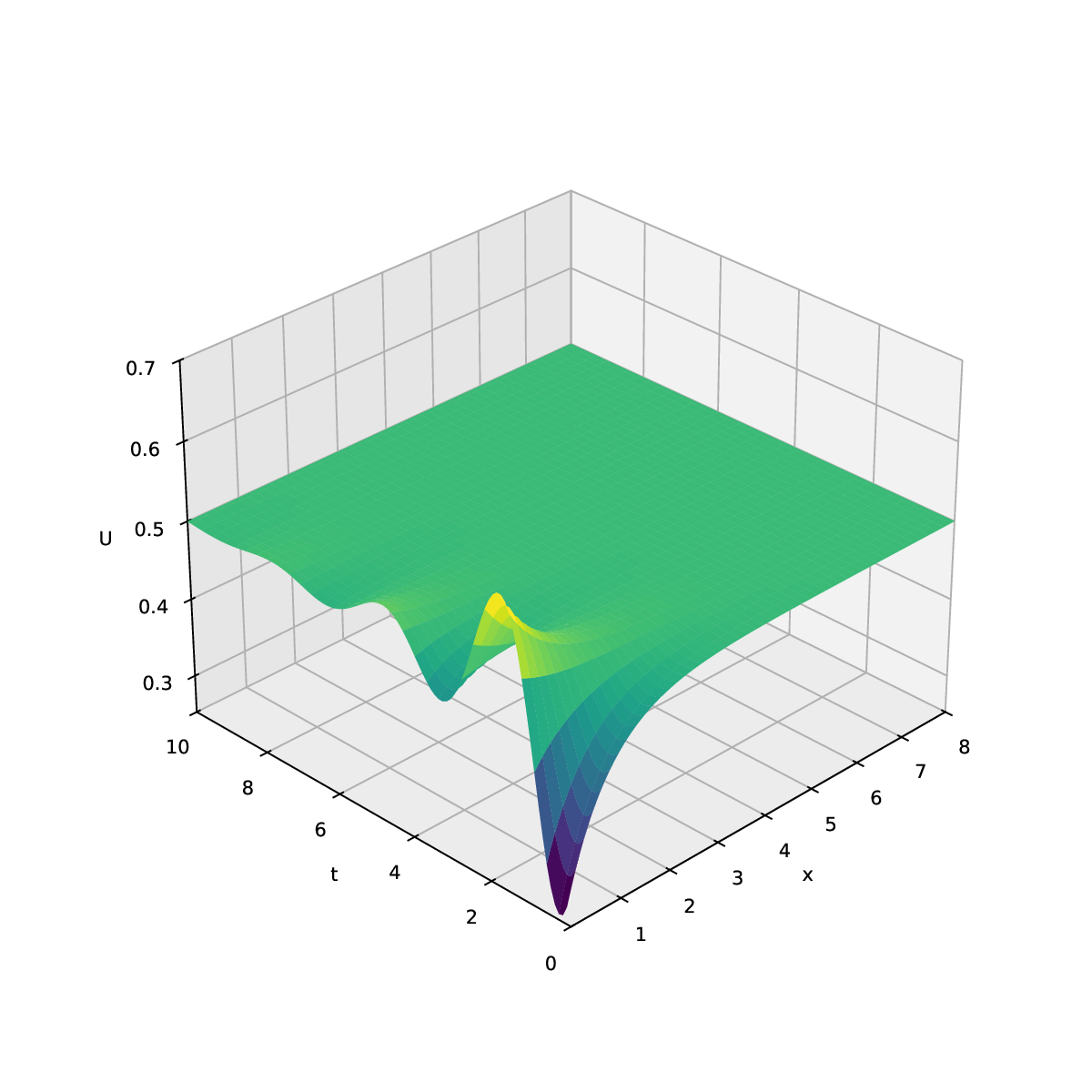}}
\end{minipage}
\hfill
\begin{minipage}[t]{9.0cm}
\centerline{\includegraphics[width=9.0cm]{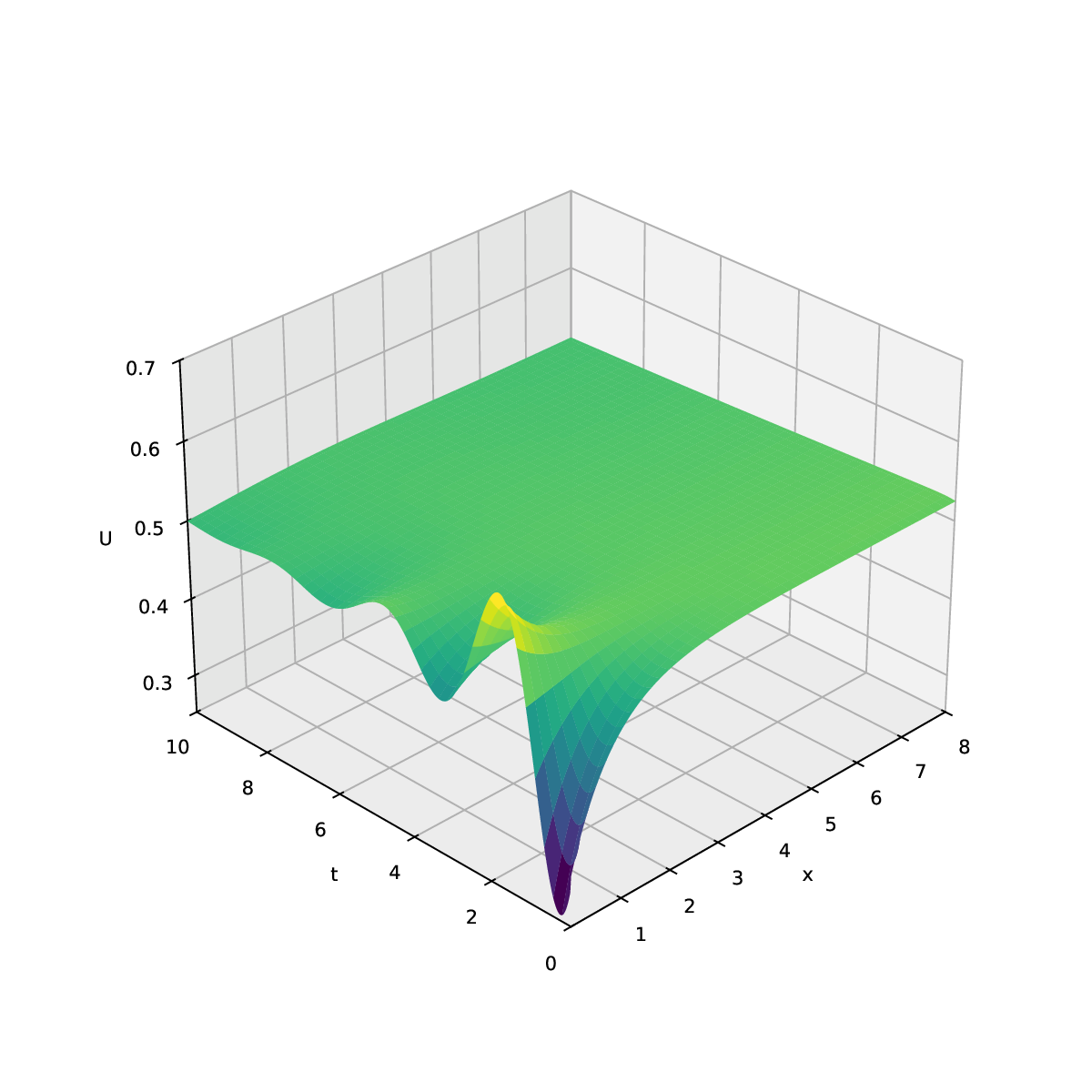}}
 \end{minipage}
\caption{Left: the component $U$ of the exact solution (\ref{5-37})
  with $e_1=e_2=0.2$. Right: the
component $U$ of the numerical solution obtained for the initial
profile (\ref{5-42}) with  $\epsilon=0.05$.} \label{Fig-4a}
 \end{figure}

\begin{figure}[t]
\begin{minipage}[t]{9.0cm}
\centerline{\includegraphics[width=9.0cm]{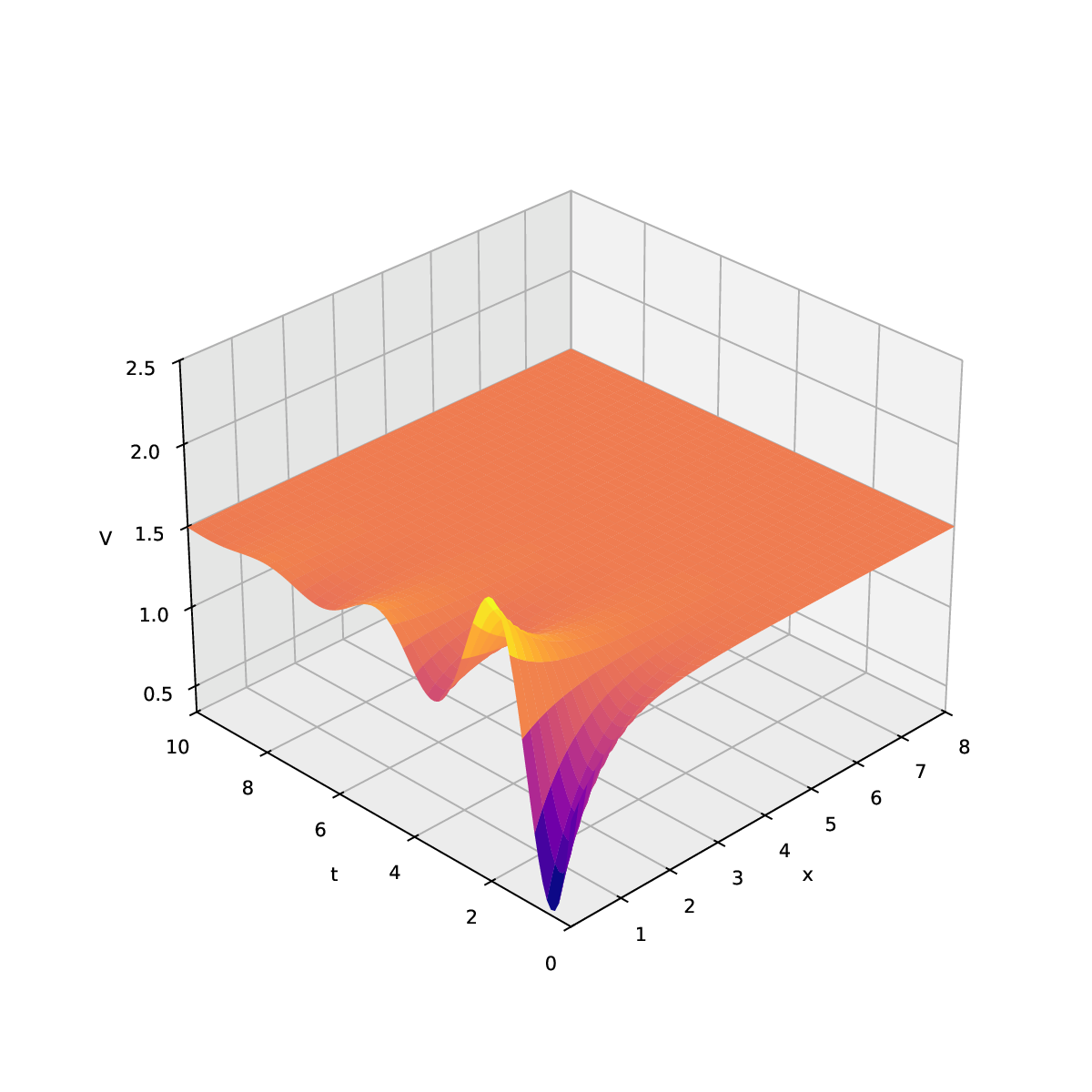}}
\end{minipage}
\hfill
\begin{minipage}[t]{9.0cm}
\centerline{\includegraphics[width=9.0cm]{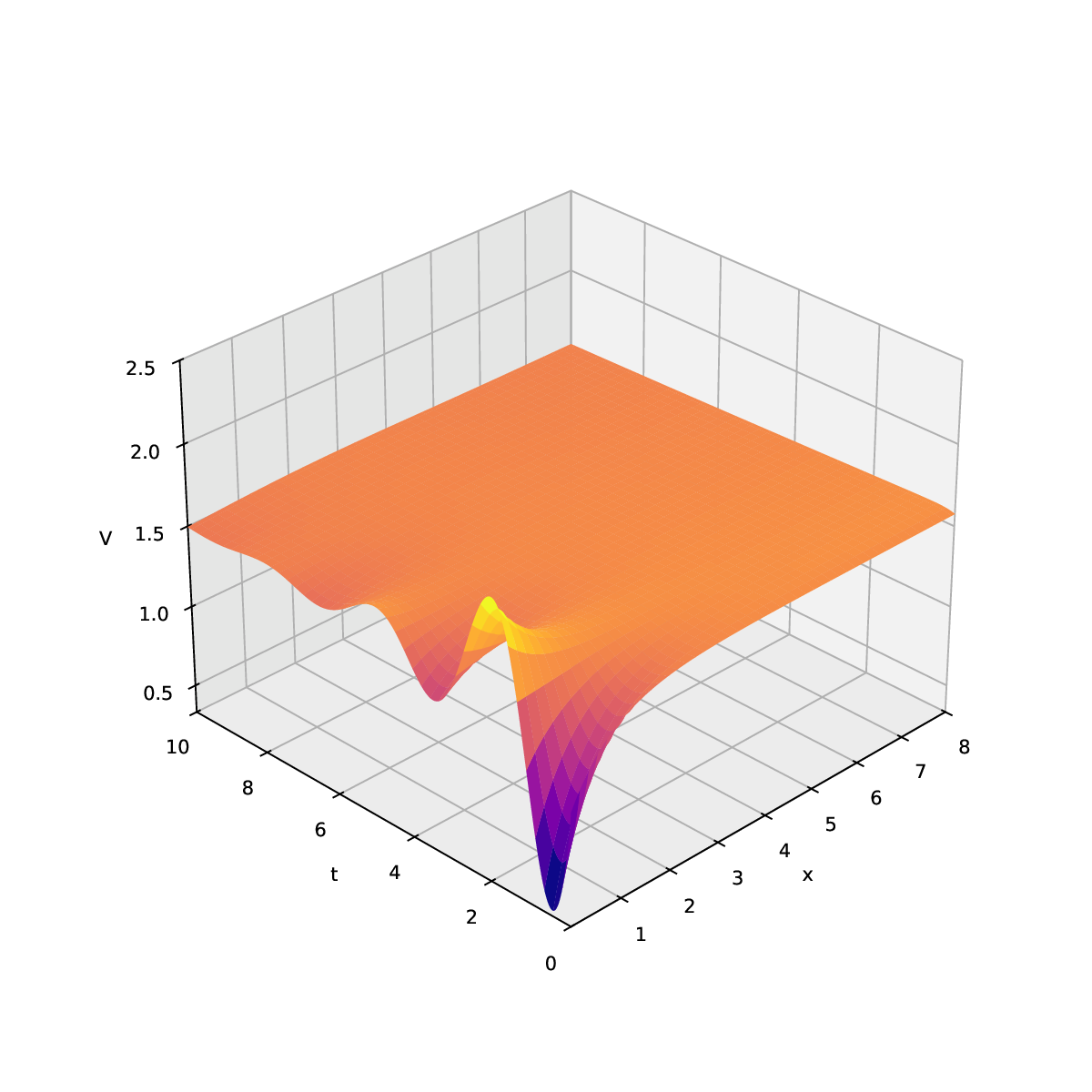}}
\end{minipage}
\caption{Left: the component $V$ of the exact solution (\ref{5-37})
  with $e_1=e_2=0.2$. Right: the
component $V$ of the numerical solution obtained for the initial
profile (\ref{5-42}) with $\epsilon=0.05$.}\label{Fig-4b}
 \end{figure}

   One observes that any  solution of the
  form (\ref{5-37}) and   (\ref{5-38}) tends
  to $(U_1, V_1)$ and  $(U_2, V_2)$, respectively, if $t \to
  +\infty$ or  $x \to +\infty$.  The exact solution (\ref{5-37})
  with $e_1=e_2=0.2$ is plotted in Fig.\ref{Fig-4a}--\ref{Fig-4b}.
If we turn to a biological interpretation
  of the exact solutions  obtained, then (\ref{5-37}) and
  (\ref{5-38})  may describe the prey-predator interaction with the
  above mentioned asymptotic behaviour and this means coexistence of
  preys and predators.


Using the properties of the exact solutions (\ref{5-37}) and
(\ref{5-38}), one may solve some boundary-value problems. Here we
present  an example in detail.

 \begin{figure}[t]
\begin{minipage}[t]{9.0cm}
\centerline{\includegraphics[width=9.0cm]{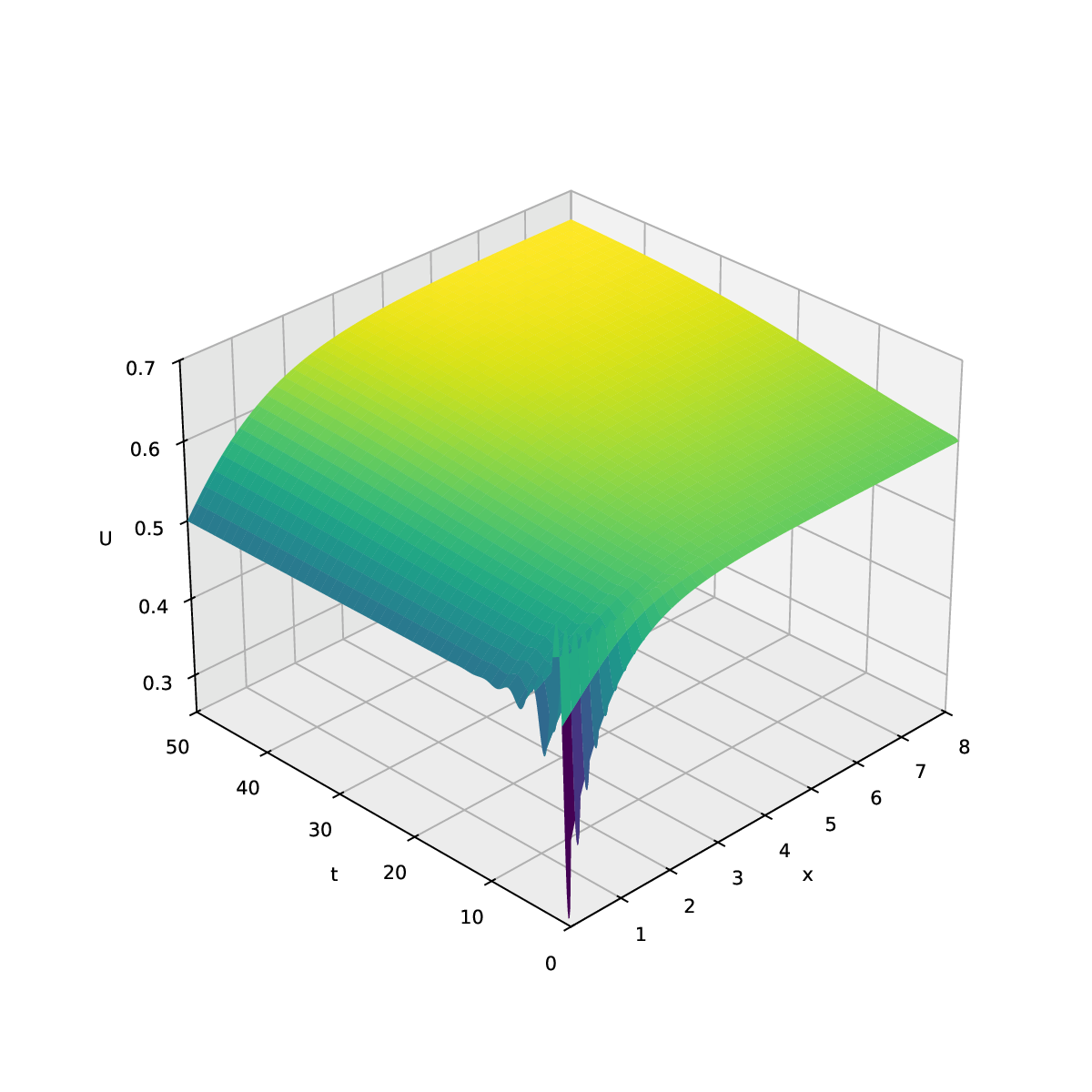}}
\end{minipage}
\hfill
\begin{minipage}[t]{7.0cm}
\centerline{\includegraphics[width=9.0cm]{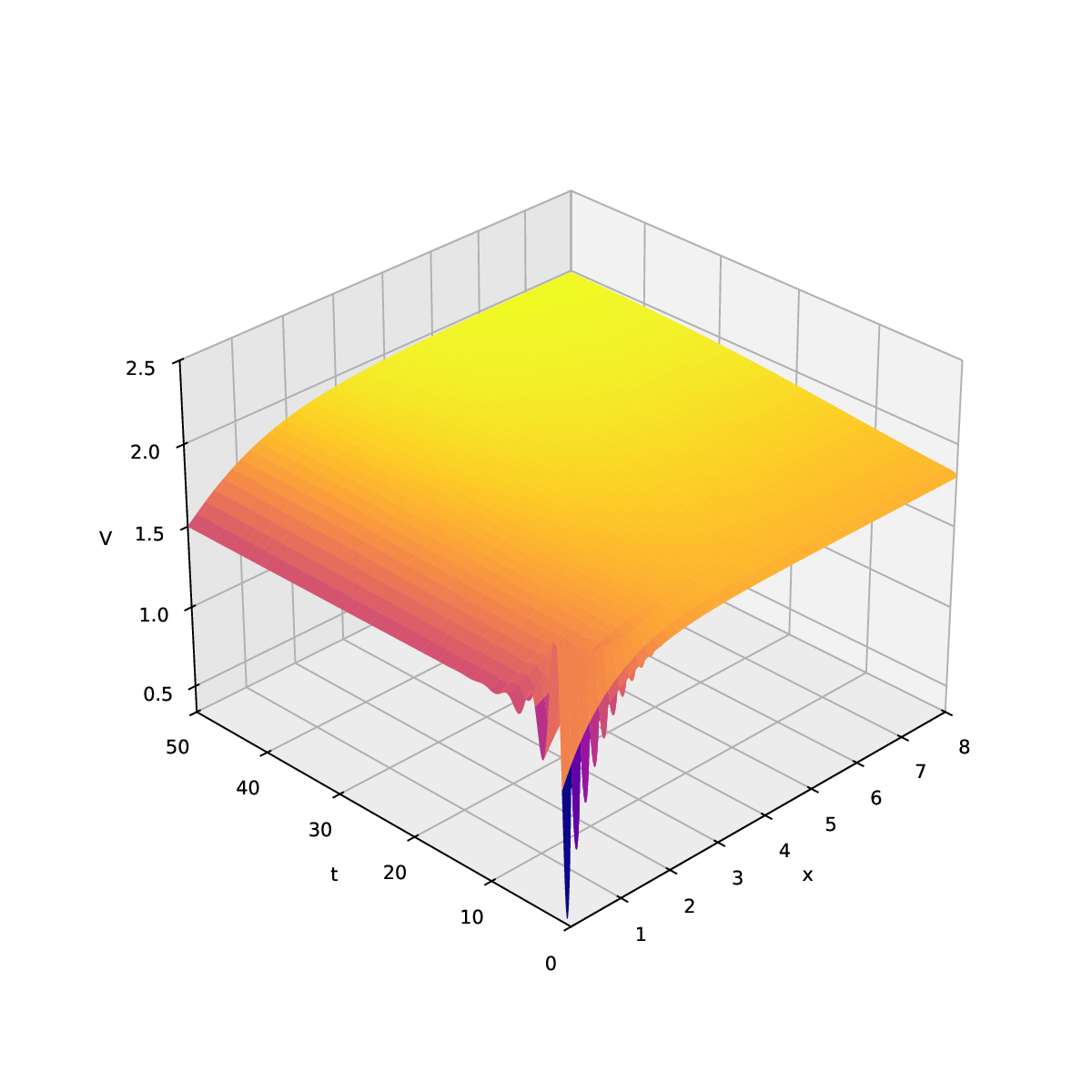}}
\end{minipage}
\caption{The numerical solution obtained for for the initial profile
(\ref{5-42}) with $\epsilon=0.20$.} \label{Fig-5b}
 \end{figure}

\bteo
 The  bounded exact solution of the
for the RD system  (\ref{5-36}) with the initial profiles
\begin{equation}\label{5-39}
 \ba
U(0,x)={1\over 2}-e_2\exp\Big(-\sqrt{\frac {3}{2}}x \Big)\equiv
U_0(x) \\
 V(0,x)={3\over 2}-\frac{1}{4}\Big(19e_2-\sqrt{71}e_1\Big)
 \exp\Big(-\sqrt{\frac {3}{2}}x\Big)\equiv  V_0(x)
 \ea
\end{equation} and the boundary   conditions
\begin{equation}\label{5-40a}
 \ba
U(t,0)={1\over 2}-\exp\Big(-{1\over 2}t \Big)\biggl(e_1\sin(\frac
{\sqrt{71}}{4}t)+e_2 \cos(\frac {\sqrt{71}}{4}t)\biggr)\\
 V(t,0)={3\over 2}-{1\over 4}\exp\Big(-{1\over 2}t
\Big)\biggl(\Big(19e_1+\sqrt{71}e_2\Big)\sin(\frac {\sqrt{71}}{4}t)\\
\qquad+\Big(19e_2-\sqrt{71}e_1\Big) \cos(\frac
{\sqrt{71}}{4}t)\biggr),
 \ea
\end{equation} and
 \be \label{5-40b}
U_x(t,+\infty)=0, \quad
  V_x(t,+\infty)=0 ,  \ee
 in the domain $ (t,x) \in (0,+ \infty )\times (0,+ \infty) $
 is given by the formulas (\ref{5-37}).\et
 {\bf Remark 4.} Setting $e_2=0$  and  $e_2=\frac
 {\sqrt{71}}{19}e_1$,
  one obtains
BVP with a constant initial profile for $U$ and  $V$, respectively.
If one sets simultaneously  $e_2=0$  and  $e_2=\frac
 {\sqrt{71}}{19}e_1$ then  formulae (\ref{5-37}) produce the steady-state
 solution.
\medskip

Now we are looking for solutions of the above BVP  with perturbed
initial conditions. We are interesting to know how  these
perturbation may affect the exact solution (\ref{5-37}).
So, let us replace (\ref{5-39}) by
 \be\label{5-42} U(0, x) = U_0(x) +\epsilon U_0(x), \quad
 V(0, x) = V_0(x) + \epsilon V_0(x), \ee  where
 $\epsilon$ is a real parameter such that $|\epsilon|<1.$

 In order to construct
numerical solutions of BVP (\ref{5-36}), (\ref{5-40a}),
(\ref{5-40b}) and (\ref{5-42}), numerical simulations were
conducted.
These simulations utilized the {\it  odeint function} from the {\it
Python  scipy.integrate package}, which internally employs {\it
LSODA} from the FORTRAN   library to solve ODEs. Notably, the
infinite space interval $(0,+ \infty )$ was replace by the finite
interval $(0,L)$ with $L=8$. Obviously, the no-flux boundary
conditions at $x=8$ for exact solutions (\ref{5-37})--(\ref{5-38})
are fulfilled with high exactness.

  Many  perturbations with different values of
  $\epsilon$
 were tested.  It turns out
  that the value of the parameter $\epsilon$ plays a crucial role on the
  solution behaviour and some  results  are  presented
  in Fig.\ref{Fig-4a}--\ref{Fig-5b}.

  Fig. \ref{Fig-4a} and \ref{Fig-4b} represent the
  components $U$ and $V$ (right-hand-side plots) of the numerical solution,
   obtained for the value
  $\epsilon=0.05$. This numerical solution has the quite similar form to
 the  exact solution  (\ref{5-37}) pictured in  Fig. \ref{Fig-4a} and \ref{Fig-4b}
 (left-hand-side plots) and practically coincides  with the latter
 excepting a very small  vicinity of $t=0$.
  Such a behaviour of  numerical solutions has been observed for any
  sufficiently small  $\epsilon$ and even for  $\epsilon \simeq 0.10$.


It turns out that the situation changes drastically in the case of
large perturbations of the   initial conditions (\ref{5-39}).
 In Fig.\ref{Fig-5b}, we present the
  components $U$ and $V$ of the numerical solution,  obtained for the value
  $\epsilon=0.20$. Nevertheless numerical solution has the similar form to
 the  exact solution  (\ref{5-37}) for small time values,
 it does not tend   to (\ref{5-37}). On the
 other hand, we observe
 that  this numerical solution has the same asymptotic behaviour as
  the exact solution (\ref{5-38}) excepting   a vicinity
  of the point $x=0$. In fact, $U \to {2\over 3}, \, V \to 2$ as
  $t \to +\infty$  for arbitrary  $ x \in [\delta, +\infty), \,
  \delta>0$.
  Obviously, there exists such a critical $\epsilon_{cr}$ that any
  smaller $\epsilon$  leads to a solution of the relevant BVP, which
  is sufficiently close to the exact solution (\ref{5-37}), while
  any larger $\epsilon$ explores a solution, which differs
  essentially from (\ref{5-37}). We are going to estimate
  analytically and numerically $\epsilon_{cr}$ in a forthcoming
  study.

  The difference between the numerical solutions plotted
   on Fig.\ref{Fig-4a}--\ref{Fig-4b}
    and Fig.\ref{Fig-5b}
  can be explained using  the $(U,V)$-phase
plane of solutions (see Fig.\ref{Fig-2b}).  In fact, the initial
profile (\ref{5-42}) with $\epsilon=0.05$
produces a point belonging to the domain $G_1$. Hence the relevant
curve of the $(U,V)$-phase plane tends to the stable node $P_1$. In
the case $\epsilon=0.2$, one observes that the relevant point moves
to the domain $G_2$, if $x$ is sufficiently large, therefore the
solution  tends to the stable node $P_2$. In the case of the
sufficiently small $x>0$, the solutions of the relevant BVP still
tend to the node $P_1$ because of the
boundary  condition (\ref{5-40a}). If this condition
 is replaced by the zero Neumann condition then
 solutions  would tend to the node $P_2$ for the arbitrary $x>0$.

  {\it Thus, the two-parameter families of exact solutions
  obtained in Section 2 play an important role for solving
   BVP (\ref{5-36}),(\ref{5-40a}) and
(\ref{5-40b})  with a wide range of initial profiles of the form
(\ref{5-39}).}

Moreover, we assume that the same situation occurs for more general
forms of initial profiles, however, relevant investigations lie
beyond the specific aims of this study.

\section{Discussion}

In this work, the  Lotka-Volterra type system with porous diffusion
(\ref{3-1}) was  studied. The system can be considered as a
simplification of the well-known SKT system (\ref{0-1})  or as an
alternative model to the classical Lotka-Volterra system with porous
diffusion. Multiparameter families of exact solutions of the system
in question
 are constructed and  their properties are established. It is shown that the solutions obtained can satisfy
  the zero Neumann conditions, which are  typical conditions for
  mathematical models describing real-world processes.

It should be noted that  plane wave solutions can be derived as
particular cases of those constructed in this work. For example,
plane wave solutions of the RD system (\ref{2-12}) can easily be
derived using formulae (\ref{4-8})  and (\ref{4-14}) by setting
either $e_1=0$, or $e_2=0$. Interestingly, the structure of the
solutions obtained is identical to that for the porous-Fisher
equation presented in \cite [Section 13.4]{mur2002}.

We also point out that the exact solutions obtained here cannot be
constructed  using the Lie method
\cite{bl-anco-10}--\cite{ch-se-pl-book},
 which is extensively applied  for finding   exact solutions of nonlinear PDEs.
 In fact, using the Lie symmetry classification of the general class
 of RD systems with nonconstant diffusivities that was  derived  in
\cite{ch-king3}, one easily checks that the RD system (\ref{3-1})
with $b_1b_2\not=0$ and  $c_1c_2\not=0$ admits only a trivial Lie
symmetry that   leads to the ansatz (\ref{2-2}).

  It is proved
  that the system  possesses two stable steady-state  points
  provided  its  coefficients  are correctly-specified. In
  particular, this occurs when the system models  the prey-predator interaction.
  A simple example of such model (see (\ref{5-2*})) is derived that
  predicts two different points of  coexistence of preys and predators, depending on the initial condition.
  It should be noted that the classical Lotka-Volterra system for
  the prey-predator interaction does not possess stable steady-state
  points.

An example of the RD system with correctly-specified coefficients
(see system (\ref{5-36})), which models the prey-predator
interaction, is studied in detail.
  The relevant  exact solutions are used for solving the system   with the mixed boundary conditions
and  initial profiles. Furthermore, a wide range of BVPs with the
above boundary conditions but perturbed initial profiles are
numerically solved.
  The exact solution  is compared with  numerical solutions.  It is  concluded   that
 the numerical  solutions  coincide
 with  the
exact solution with high exactness provided perturbations  of the
initial profiles are sufficiently small. In the case of large
perturbations  of the initial profiles, the relevant numerical
solutions differ essentially from the exact solution. In particular,
these solutions possess  different asymptotical behaviour. Thus, the
exact solutions obtained  play an important role for solving some
boundary-value problems for  the  RD system (\ref{3-1}). Finally, we
note  that a similar  investigation for  the classical
Lotka-Volterra system with linear diffusion was performed in
\cite{ch-du-2004}.

\section{Acknowledgements}
  R.Ch. acknowledges
that this research was  funded by  the British Academy (Leverhulme
Researchers at Risk Research Support Grant LTRSF-24-100101). The
authors are grateful to Dr Vasyl' Dutka for fruitful discussions
concerning numerical simulations.


\end{document}